\documentclass{article}
\begin{document}

\title{Minimizing Squared Perpendicular Errors}
\author{Donald E. Ramirez \and der@virginia.edu \and University of Virginia
\and Department of Mathematics \and P. O. Box 400137 \and Charlottesville,
VA 22904-4137}
\date{}
\maketitle

\begin{abstract}
The best fit line from minimizing the sum of the squared perpendicular error
is given.
\end{abstract}

\section{Introduction}

With simple linear regression and data $\{(x_{1},y_{1}),...,(x_{n},y_{n})\},$
we will minimize the sum of the squared vertical errors to find the
"best-fit" line.

\section{Minimizing Squared Perpendicular Errors}

For a right triangle with legs $\{a,b\}$ and hypotenuse $c$ with $h$ the
distance from the intersection of the legs to the hypotenuse. From similar
triangles $a/h=c/b,$ so $h^{2}=(ab)^{2}/(a^{2}+b^{2}).$ The sum of the
squared perpendicular errors $SSE_{p}$ for the fitted line $h(x)=\beta
_{0}+\beta _{1}x$ has%
\begin{equation}
SSE_{p}(\beta _{0},\beta _{1})=\sum_{i=1}^{n}\frac{(y_{i}-\beta _{0}-\beta
_{1}x_{i})^{2}(x_{i}-(y_{i}-\beta _{0})/\beta _{1})^{2}}{(y_{i}-\beta
_{0}-\beta _{1}x_{i})^{2}+(x_{i}-(y_{i}-\beta _{0})/\beta _{1})^{2}}
\label{SSEp}
\end{equation}%
with 
\begin{equation}
\frac{\partial }{\partial \beta _{0}}SSE_{p}(\beta _{0},\beta _{1})=\frac{2}{%
\beta _{1}^{2}+1}\left( \beta _{1}\sum_{i=1}^{n}x_{i}+n\beta
_{0}-\sum_{i=1}^{n}y_{i}\right)
\end{equation}%
with root \ 
\begin{equation}
\widehat{\beta }_{0}=\overline{y}-\widehat{\beta }_{1}\overline{x}
\label{b0 for perpendicular errors}
\end{equation}%
as in the simple linear regression case

Set $S_{xx}=\sum_{i=1}^{n}(x_{i}-\overline{x})^{2}$, $S_{yy}=%
\sum_{i=1}^{n}(y_{i}-\overline{y})^{2}$, $S_{xy}=\sum_{i=1}^{n}(x_{i}-%
\overline{x})(y_{i}-\overline{y})$. With $\beta _{0}$ from Equation \ref{b0
for perpendicular errors}, 
\begin{equation}
SSE_{p}(\beta _{1})=\frac{S_{yy}-2\beta _{1}S_{xy}+\beta _{1}^{2}S_{xx}}{%
1+\beta _{1}^{2}},  \label{SSE}
\end{equation}%
and the equation for the slope $\beta _{1}$ is given by 
\begin{equation}
\frac{\partial }{\partial \beta _{1}}SSE_{p}(\beta _{1})=\frac{2\left( \beta
_{1}^{2}S_{xy}+\beta _{1}(S_{xx}-S_{yy}\right) -S_{xy})}{(1+\beta
_{1}^{2})^{2}}  \label{dSSE/db1}
\end{equation}%
with roots (provided $S_{xy}\neq 0)$ 
\begin{equation}
\widehat{\beta }_{1}=\frac{S_{yy}-S_{xx}\pm \sqrt{%
(S_{yy}-S_{xx})^{2}+4S_{xy}^{2}}}{2S_{xy}}.
\label{b1 from quadratic formula}
\end{equation}

If $S_{xy}\neq 0,$ then then there are two roots from Equation \ref{b1 from
quadratic formula}, one positive and one negative. If $S_{xy}>0$, then the
minimum of $SSE_{p}(\beta _{1})$ is the positive root, and if $S_{xy}<0$,
the minimum is the negative root. This follows from the fact that there are
two critical points for $\frac{\partial }{\partial \beta _{1}}SSE_{p}(\beta
_{1})$ one to the left of zero and one to the right of zero, and that the
value at zero of $\frac{\partial }{\partial \beta _{1}}SSE_{p}(0)=-2S_{xy}.$

If $S_{xy}=0,$ the minimum of $SSE_{p}(\beta _{1})$ is found from Equation %
\ref{SSE} and (1) is $S_{yy}$ with $\beta _{1}=0$ for $S_{yy}<S_{xx};$ (2)
is $S_{xx}$ with $\beta _{1}=\infty $ for $S_{xx}<S_{yy};$ and (3) is $%
S_{xx}=S_{yy}$ with any value of $\beta _{1}$ for $S_{xx}=S_{yy}.$

\section{Example}

For the data $\{(0,0),(1,1),(1,0),(0,0)\}$ with $\{\overline{x}=1/2,%
\overline{y}=1/4,S_{xx}=1,S_{yy}=3/4,S_{xy}=1/2,\rho =\sqrt{3}/3\}$, the
"best" fit straight line has $\widehat{\beta }_{1}=0.78078$ and $\widehat{%
\beta }_{0}=-0.14039$.

\section{Summary}

We have investigated the fitted linear equation to data using as the
criterion for fit the minimized sum of squares of the perpendicular errors.

\end{document}